\documentclass[12pt]{amsart}
\usepackage{amssymb,amsmath}
\numberwithin{equation}{section}
\def\simleq{\underset\sim<}
\def\simgeq{\underset\sim>}
\def\simle{\underset\sim<}
\def\simge{\underset\sim>}
\def\T{\text}
\def\1#1{\overline{#1}}
\def\2#1{\widetilde{#1}}
\def\3#1{\widehat{#1}}
\def\4#1{\mathbb{#1}}
\def\5#1{\frak{#1}}
\def\6#1{{\mathcal{#1}}}
\def\C{{\4C}}
\def\R{{\4R}}

\def\sumK{\underset{|K|=k-1}{{\sum}'}}
\def\sumJ{\underset{|J|=k}{{\sum}'}}

\def\Re{{\sf Re}\,}

\def\phi{\varphi}

\newtheorem{Thm}{Theorem}[section]
\newtheorem{Cor}[Thm]{Corollary}
\newtheorem{Pro}[Thm]{Proposition}
\newtheorem{Lem}[Thm]{Lemma}
\theoremstyle{definition}\newtheorem{Def}[Thm]{Definition}
\theoremstyle{remark}
\newtheorem{Rem}[Thm]{Remark}
\newtheorem{Exa}[Thm]{Example}

\def\Label#1{\label{#1}}
\def\bl{\begin{Lem}}
\def\el{\end{Lem}}
\def\bp{\begin{Pro}}
\def\ep{\end{Pro}}
\def\bt{\begin{Thm}}
\def\et{\end{Thm}}
\def\bc{\begin{Cor}}
\def\ec{\end{Cor}}
\def\bd{\begin{Def}}
\def\ed{\end{Def}}
\def\br{\begin{Rem}}
\def\er{\end{Rem}}
\def\be{\begin{Exa}}
\def\ee{\end{Exa}}
\def\bpf{\begin{proof}}
\def\epf{\end{proof}}
\def\ben{\begin{enumerate}}
\def\een{\end{enumerate}}

\def\1alpha{[\frac1\alpha]}
\def\T{\text}
\def\R{{\Bbb R}}

\def\C{{\Bbb C}}

\numberwithin{equation}{section}
\def\T{\text}
\newcommand{\om}{\omega}
\newcommand{\bom}{\bar{\omega}}
\newcommand{\Dom}{\text{Dom} }
\newcommand{\we}{\wedge}
\newcommand{\no}[1]{\|{#1}\|}

\newtheorem{theorem}{Theorem  }[section]
\newtheorem{definition}[theorem]{Definition }
\newtheorem{lemma}[theorem]{Lemma  }
\newtheorem{proposition}[theorem]{Proposition  }
\newtheorem{corollary}[theorem]{Corollary }
\newtheorem{example}[theorem]{\it Example }

\begin{document}
\title[Subelliptiity of a weakly $Q$-pseudoconvex/concave domain]
{Subellipticity of the $\bar{\partial}$-Neumann problem on a weakly $Q$-pseudoconvex/concave domain}
\author[ T.V.~Khanh and G.~Zampieri ]
{Tran Vu Khanh and Giuseppe Zampieri}
\address{Dipartimento di Matematica, Universit\`a di Padova, via 
Trieste 63, 35121 Padova, Italy}
\email{khanh@math.unipd.it, 
zampieri@math.unipd.it}
\maketitle
\begin{abstract}
For a domain $D$ of $\mathbb{C}^n$ which is weakly $q$-pseudoconvex or $q$-pseudoconcave we   give a sufficient condition for subelliptic estimates for the $\bar{\partial}$-Neumann problem. The paper extends to domains which are not necessarily pseudoconvex, the results and the techniques of Catlin \cite{C87}.
\end{abstract}

\section{Introduction}
Let $D$ be a bounded domain of $\C^n$ with smooth boundary. For a form $f $ of degree $k$ which satisfies $\bar\partial f=0$, to solve the $\bar\partial$-Neumann  problem consists in finding a form of degree $k-1$ such that
\begin{equation}
\Label{1.1}
\begin{cases}
\bar\partial u=f,
\\
f\T{ is orthogonal to $\T{Ker}\, \bar\partial$}.
\end{cases}
\end{equation}
The main interest relies in the regularity at the boundary for this problem, that is, in stating under which condition $u$ inherits from $f$ the smoothness at the boundary $\partial D$ (it certainly does in the interior).
Let $\bar\partial^*$ be the formal adjoint of $\bar\partial$ under the choice of a smoothly varying hermitian metric on $\bar D$. Related to \eqref{1.1} is the problem  

\begin{eqnarray}
\Label{1.2}
\begin{cases}(\bar{\partial}\bar{\partial}^*+\bar{\partial}^*\bar{\partial})u=f  
\\
u\in D_{\bar{\partial}}\cap D_{\bar{\partial}^*}
\\
\bar{\partial} u\in D_{\bar{\partial}^*}, \bar{\partial}^* u\in D_{\bar{\partial}}, 
\end{cases}
\end{eqnarray}
where $D_{\bar\partial^*}$ and $D_{\bar\partial}$ are the domains of $\bar\partial^*$ and $\bar\partial$ respectively.
This is a non-elliptic boundary value problem; in fact,   the Kohn Laplacian  $\Box=\bar\partial\bar\partial ^*+\bar\partial^*\bar\partial$ itself is elliptic but the boundary conditions which are imposed by the membership to $D_\Box$ are not. If \eqref{1.1} has a solution for every $f$, then one defines the $\bar\partial$-Neumann operator $N:=\Box^{-1}$; this commutes both to $\bar\partial$ and $\bar\partial^*$. If we then turn back to \eqref{1.1} and define $u:=\bar\partial^*Nf$ we see that
\begin{equation*}
\begin{split}
\bar\partial u&=\bar\partial^*Nf
\\
&=\Box Nf=f.
\end{split}
\end{equation*}
Also, $\bar\partial^*u=\bar\partial^*\bar\partial^*Nf=0$ and therefore $u$ is orthogonal to $\T{Ker}\, \bar\partial$. One of the main methods used in investigating the regularity at the boundary of the solutions of \eqref{1.1} consists in certain a priori subelliptic estimates.
\bd
 The $\bar{\partial}$-Neumann problem is said to satisfy a subelliptic estimate of order $\epsilon>0$ at $z_o\in \bar D$ on $k$ forms if there exist a positive constant $c$ and a neighborhood $V\ni z_o$ such that 
\begin{eqnarray}
\Label{1.3,5}
||u||_\epsilon \le c(\no{\bar{\partial} u}^2+\no{\bar{\partial}^* u}^2)\quad\T{ for any }u\in C^\infty_c(\bar{D}\cap V)^{k}\cap D_{\bar{\partial}^*}.
\end{eqnarray}
\ed
By  Garding's inequality,  subelliptic estimates of order $1$, that is, elliptic estimates  hold in the interior of $D$. So our interest is confined to the boundary $\partial D$. When the domain $D$ is pseudoconvex, a great deal of work has been done about subelliptic estimates. The most general results concerning this problem have been  obtained by Kohn \cite{K79} and  Catlin \cite{C87}.  

In \cite{K79}, Kohn gave a sufficient condition for subellipticity over pseudoconvex domains with real analytic boundary       by introducing a sequence of ideals of subelliptic multipliers.

In \cite{C87}, Catlin proved, regardless whether $\partial D$ is real analytic or not, that subelliptic estimates hold  for $k$ forms at  $z_o$ if and only if a certain number $D_k(z_o)$ is finite. Note that the definition of $D_k(z_o)$ in \cite{C87} is closely related to that of $\Delta_k(z_o)$ due to   D'Angelo. In particular, when $k=1$, these numbers do coincide.

However, not much is known in the case when the domain is not necessarily pseudoconvex except from the results
related to the celebrated $Z(k)$ condition which characterizes the existence of subelliptic estimates for $\epsilon=\frac12$ 
according to H\"ormander \cite{H65} and Folland-Kohn \cite{FK72}. Some further results, mainly related to the case of forms of top degree $n-1$ are due to Ho \cite{Ho85}.

We exploit here the full strength of Catlin's method to study  subellipticity on domains which are not necessarily pseudoconvex.
Let $\partial D$ be defined by $r=0$ with $r<0$ on the side of  $D$ and let $T^\C\partial D$ be the complex tangent bundle to $\partial D$.  We use the following notations: $L_{\partial D}=(r_{ij})|_{T^\C\partial D}$  is the Levi form of the boundary,  $s^+_{\partial D },\,\,s^-_{\partial D},\,\,s^0_{\partial D}$  are the numbers of eigenvalues of $L_{\partial D}$ which are $>0,\,\,<0,\,\,-0$ respectively and finally  $\lambda_1\leq \lambda_2\leq...\lambda_{n-1}$ are its  ordered eigenvalues.
We take a pair of indices   $1\le q\le n-1$ and $0\le q_o\le n-1$ such that $q\neq q_o$.
We assume that there is a bundle $\mathcal V^{q_o} \subset T^{1,0}\partial D$ of  rank $q_o$ with smooth coefficients in $V$, say the bundle of the first $q_o$ coordinate tangential vector fields $L_1,...L_{q_o}$, such that 
\begin{equation}
\Label{1.3}
\underset{j=1}{\overset q\sum}\lambda_j-\underset {j=1}{\overset {q_o}\sum}r_{jj}\ge 0\quad\T{ on $\partial D$}.
\end{equation}
\bd
(i) If $q>q_o$ we say that $D$ is $q$-pseudoconvex.

\noindent
(ii) If $q<q_o$ we say that $D$ is $q$-pseudoconcave.
\ed
\br
The notion of $q$-pseudoconvexity was used in \cite{Z00}  to prove the existence of $C^\infty(\bar D)$ solutions to the equation $\bar\partial u=f$.  Though the notion of $q$-pseudoconcavity is formally simmetric to $q$-pseudoconvexity, it is useless in the existence problem. The reason is intrinsic. Existence is a  ``global" problem
but  bounded domains are never globally $q$-pseudoconcave.
Owing to the local nature of subelliptic estimates and the related hypoellipticity of $\bar\partial$, this is the first occurence where $q$-pseudoconcavity comes successfully into play.
\er
\br
Assume that \eqref{1.3} holds. Then, if $q>q_o$, we must have $\lambda_q\ge 0$.  Thus \eqref{1.3} still holds with the first sum $\underset{j=1}{\overset q\sum}$ replaced by $\underset{j=1}{\overset k\sum}$ for any $k\ge q$.
Similarly, if $q<q_o$, then $\lambda_{q+1}\le 0$. Thus  \eqref{1.3}  holds with  $\underset{j=1}{\overset q\sum}$ replaced by $\underset{j=1}{\overset k\sum}$ for any $k\le q$.
\er
\br
When we have strict inequality ``$<$" in \eqref{1.3}, it means that we have in fact $\lambda_q>0$ and $\lambda_{q+1}<0$ in the two cases respective $q>q_o$ and $q<q_o$. It follows
\begin{equation}
\Label{1.4}
q>n-1-s^+\quad(\T{resp. $q<s^-$}).
\end{equation}
We refer to these two situations as ``strong" $q$-pseudoconvexity (resp. -pseudoconcavity). Note that this is the same as to saying, in the terminology of Folland-Kohn,  that $D$ satisfies $Z(k)$ for any $k\ge q$ (resp. $k\le q$).
\er
We write $k$-forms as $u=(u_J)_J$ where $J=j_1<j_2<...<j_k$ are ordered multiindices. When the multiindices are not ordered, the coefficients are assumed to be alternant. Thus, if $J$ decomposes as $J=jK$, then $u_{jK}=\T{sign}\binom{J}{jK}u_J$.
We define the $\delta$-strip of $D$ along the boundary   by $S_\delta=\{z\in D: r(z)>-\delta\}$.
The main result in this paper is the following.
\bt
\Label{t1.1}
Let \eqref{1.3} be satisfied in a neighborhood of $z_o$,
let $k\geq q$ (resp. $k\leq q$) for $q> q_o$ (resp. $q<q_o$)
 and suppose that for small $\delta$ there exists a weight $\phi=\phi^\delta$ in $C^2(\bar D\cap V)$ such that
\begin{equation}
\Label{1.5}
\begin{cases}\phi\le 1, 
\\
\begin{split}
{\sum_{|K|=k-1}}'\sum^n_{i,j=1}\phi_{ij}(z)u_{iK}\bar{u}_{jK}&-{\sum_{|J|=k}}'\sum^{q_o}_{j=1}\phi_{jj}|u_{J}|^2
\\
&\ge c\underset{j=1}{\overset{q_o}\sum}\left| L_j(\phi)\right|^2|u|^2,
\T{ for any } z\in \bar{D}\cap V
\end{split}
\end{cases}
\end{equation}
and
\begin{equation}
\Label{1.6}
\begin{cases}|\phi|\le 1, 
\\
\begin{split}
{\sum_{|K|=k-1}}'\sum^n_{i,j=1}\phi_{ij}(z)u_{iK}\bar{u}_{jK}&-{\sum_{|J|=k}}'\sum^{q_o}_{j=1}\phi_{jj}|u_{J}|^2
\\
&
\ge c\delta^{-2\epsilon} |u|^2
\T{ for any } z\in \bar{S_\delta}\cap V,
\end{split}
\end{cases}
\end{equation}
where the constant $c>0$ does not depend on $\delta$ or $u$.

Then, $\epsilon$-subelliptic estimates at $z_o$ hold for forms of degree $k$.
\et
%Note that if $M$ is strongly $q$-pseudoconvex (respectively -pseudoconcave) that is if $q>n-1-s^+$ (resp. $q<s^-$), then $Z(k)$ holds for $k\ge q$ (resp. $k\le q$). In this case, as it has already been remarked, $\frac12$-subelliptic estimates in corresponding degrees $k$ are classical.
%We now show how to get the weight $\phi^\delta$ which occurs in Theorem~\ref{t1.1}. As an application, we get subelliptic estimates in degree $k\ge q$ (resp. $k\le q$) which generalize the conclusions of Ho \cite{Ho85} which are confind to the case $k=n-1$.
It is not restrictive to assume, as we will do all throughout the paper, that $L_j(z_o)=\partial_{z_j}$ for any $j$.
For every $q$-pseudoconvex/concave domain there is a small perturbation for which subelliptic estimates hold.
\bt
\Label{t1.2}
Let $D$ be $q$-pseudoconvex (resp. $q$-pseudoconcave); thus it is defined by $r<0$ for $r:=2\Re z_n+h(z_1,...,z_{n-1},y_n)$ satisfying \eqref{1.3} for $q>q_o$ (resp. $q<q_o$). Let $\tilde r:=r+\sum^{n-1}_{j=q} h_j(z_j)$ (resp. $\tilde r=r-\underset{j=1}{\overset {q+1} \sum }h_j(z_j)$) where the $h_j$'s are 
real positive subharmonic, non harmonic, functions 
 of vanishing order $2m_j$ that, by reordering, we may assume to be decreasing $...m_j\ge m_{j+1}...$ (resp. increasing $...m_j\le m_{j+1}...$) and let $\tilde D$ be defined by $\tilde r<0$.

Then subelliptic estimates hold for $\tilde D$ in degree $k\ge q$ (resp. $k\le q$) of any order $<\epsilon_k$ for $\epsilon_k:=\frac1{2m_{k}}$  (resp. $\epsilon_k:=\frac1{2m_{k+1}}$).
In both cases, when $\epsilon_k=\frac12$, we have in fact estimates including for order $\frac12$.
\et
When $\epsilon_k=\frac12$ it means that $Z(k)$ is satisfied; thus we regain the result by H\"ormander and Folland-Kohn.
We will refer to functions such as the above $h_j$'s as subharmonic functions satisfying $h_j\cong |z_j|^{2m_j}$.
\be
Let $D$ be defined by
$$
2\Re z_n-\underset{j=1}{\overset {q_o}\sum}|z_j|^{2m_j}+\underset{j=q_o+1}{\overset{n-1}\sum}|z_j|^{2m_j}<0,
$$
where the two groups of indices $\{m_1,...,m_{q_o}\}$ and $\{m_{q_o+1},...,m_{n-1}\}$ 
have increasing and decreasing order respectively.
Then subelliptic estimates hold in  degree $k\neq q_o$ of any order smaller than $\epsilon_k$  defined in 
Theorem~\ref{t1.2}. 
\ee

\bc
\Label{c1.1}
Let $D$ be a domain in $\mathbb{C}^n$ defined by 
$$2\Re z_n+g+|z_{n-1}|^{2m}<0$$
where $g$ is a real $C^\infty$ function such that $g_{n-1\,n-1}=o(|z_{n-1}|^{2(m-1)})$. Then subelliptic estimates of order $\epsilon<\frac{1}{2m}$ hold at $z_o=0$ for any $(n-1)$-form.
\ec

\bpf
Put $r:=2\Re z_n+g+\frac12|z_{n-1}|^{2m}$; we claim that  $r$ satisfies \eqref{1.3} for $q_o=n-2$ and $q=n-1$.
In fact for a tangential $(n-1)$-from the only coefficient which does not vanish at $\partial D$ is $u_J=u_{1,...,n-1}$. Thus
$$
\sumK\underset{ij=1}{\overset {n-1}\sum}r_{ij}u_{iK}\bar u_{jK}-\sumJ\underset{j=1}{\overset {n-2}\sum}r_{jj}|u_J|^2=(|z_{n-1}|^{2(m-1)}+o(|z_{n-1}|^{2(m-1)})|u_{1,...,n-1}|^2
$$
 which is $\ge0$.
We are thus in position to
apply Theorem~\ref{t1.2}.

\epf
\be
Let $D$ be defined by 
$$2\Re z_3-|z_1^2+z_2^3|^{2} \pm|z_1|^{2m}+|z_2|^4 <0\,\,\text{ or  }\,\, 2\Re z_3-|z^2_1z_2^3|^{2}\pm|z_1|^{2m}+|z_2|^4<0;$$
then subelliptic estimates  hold at $z_o=0$ on 2-forms for any order $\epsilon<\frac14$.
\ee

Remark : Corollary~\ref{c1.1} is more general than Corollary 3.4 in \cite{Ho91} where $g$ cannot depend on $z_{n-1}$ and $y_n$.

 We decompose the coordinates as $z=(z',z'',z_n)\in\C^{q_o}\times\C^{n-q_o-1}\times\C$. The conclusion contained in  
 Theorem~\ref{t1.2} is sharp.
\bt
\Label{t1.3}
(i) Let $r=2\Re z_n-Q(z')$ for $Q\ge0$ and set $\tilde r=r+\underset{j=q_o+1}{\overset{n-1}\sum}h_j(z_j)$ where the $h_j$'s are subharmonic and satisfy $h_j\cong |z_j|^{2m_j}$ with $m_j\ge m_{j+1}\ge...$ (decreasing)  and with $Q=O(|z'|^{2m_{q_o+1}})$. If $\epsilon$-subelliptic estimates at $z_o=0$ hold in degree $k\ge q_o+1$, then we must have $\epsilon\le \frac1{2m_k}$. 

\noindent
(ii) Let $r=2\Re z_n+Q(z'')$ for $Q\ge0$ and set $\tilde r=r+\underset{j=1}{\overset {q_o}\sum}h_j(z_j)$ with $h_j$ subharmonic satisfying $h_j\cong |z_j|^{2m_j}$ with $m_j\le m_{j+1}\le...$ (increasing). We also assume $m_1\ge \frac{m_{q_o-1}}2+\frac14$ and $Q=O(|z''|^{2m_{q_o-1}})$. 
If $\epsilon$-subelliptic estimates hold at $z_o=0$ in degree $k\le q_o-1$, then $\epsilon\le \frac1{2m_k}$.
\et
Necessary conditions for subellipticity in degree $k=n-1$ are also stated in \cite{Ho85}; however, the $\bar\partial$-Neumann conditions seem not to be respected in the proof therein.

The paper is structured as follows. In Section 2 we introduce the geometric concept of $q$-pseudoconvexity and $q$-pseudoconcavity. In Section 3 we derive some basic inequalities  which are  useful for the proof of Theorem 1.2. Sections 4, 5 and 6 are devoted to the  proof of Theorem~\ref{t1.1}, Theorem~\ref{t1.2} and Theorem~\ref{t1.3} respectively.

\section{q-pseudoconvex/pseudoconcave domains}

Let $D$ be a bounded domain in $\mathbb{C}^n$ with smooth boundary $\partial D$ defined by $r=0$ with $\partial r\neq0$. For a given boundary point $z_o\in \partial D$, we consider a  complex frame adapted to $\partial D$, that is, an orthonormal basis $\om_1,...,\om_n=\partial r$ of  $(1,0)$ forms with $C^\infty$ coefficients in a neighborhood of $z_o$. We denote by $(r_{jk}(z))^{n}_{j,k=1}$ the matrix of the Levi form $\partial\bar{\partial} r(z)$  with respect to the basis $\om^1,...,\om^n$. Let $\lambda_1(z)\le ...\le \lambda_{n-1}$ be the eigenvalues of $(r_{jk}(z))^{n-1}_{j,k=1}$ and denote $s^+_{\partial D},s^-_{\partial D},s^0_{\partial D} $ their number according to the different sign. 
Let $q$ and $q_o$ be a pair of indices for which \eqref{1.3} is fulfilled in a suitable choice of the frame; remember that we have defined $D$ to be $q$-pseudoconvex or $q$-pseudoconcave according to $q>q_o$ or $q<q_o$; for the case $q>q_o$ this definition 
follows  \cite{Z00}. The pseudoconvexity/concavity is said to be strong when \eqref{1.3} holds as strict inequality. 

As it has already been noticed, \eqref{1.3} for $q>q_o$  implies $\lambda_{q}\ge 0$; hence \eqref{1.3} is still true if we replace the first sum $\sum^{q}_{j=1}\cdot$ by $\sum^{k}_{j=1}\cdot$ for any $k$ such that $q\le k\le n-1$.
Similarly, if it holds for $q<q_o$, then $\lambda_{q+1}\le0$ and hence it also holds with $q$ replaced by $k\le q$ in the first sum.
\be
\Label{e2.1}
It is readily seen that for $q_o=s^-+s^0$ and for any $q>q_o$ (resp. $q_o=s^-$ and any $q<q_o$), \eqref{1.3} is satisfied in a suitable local boundary frame. 
Thus any index $q\notin [s^-,\,s^-+s^0]$ satisfies \eqref{1.3} for either choice of $q_o$. The interesting point in \eqref{1.3} is to be capable to get \eqref{1.3} for indices $q\notin [s^-,\,s^-+s^0]$. 
\ee
\be
Let $s^-(z)$ be constant for $z\in\partial D$ close to $z_o$; then \eqref{1.3} holds for $q_o=s^-$ and $q=s^-+1$. In fact, we have $\lambda_{s^-}<0\le \lambda_{s^-+1}$, and therefore the negative eigenvectors span a bunlde $\mathcal{V}^{q_o}$ for $q_o=s^-$ that, identified with the span of the first $q_o$ coordinate vector fields, yields $\sum^{q_o+1}_{j=1}\lambda_j(z)\ge \sum^{q_o}_{j=1} r_{jj}(z)$. Note that a pseudoconvex domain is characterized by $s^-(z)\equiv0$, thus, it is 1-pseudoconvex in our terminology. 

In the same way, if $s^+(z)$ is constant at $z_o$, then $\lambda_{s^-+s^0}\le0<\lambda_{s^-+s^0+1}$. Then, the eigenspace of the eigenvectors $\le0$ is a bundle which, identified to that of the first $q_o=s^-+s^0$ coordinate vector fields yields \eqref{1.3} for $q=q_o-1$. In particular a pseudoconcave domain, that is a domain which satisfies $s^+\equiv0$, is $(n-2)$-pseudoconcave in our terminology.
\ee

The following lemma plays an essential role in the following.
\bl
\Label{l2.1}
Let $D$ be a smoothly bounded domain. Then $D$ is $q$-pseudoconvex (resp. $q$-pseudoconcave) if and only if
\begin{eqnarray}
{\sum_{|K|=k-1}}'\sum^{n-1}_{i,j=1}r_{ij} u_{iK}\bar{u}_{jK}-{\sum_{|J|=k}}'\sum^{q_o}_{i=1}r_{jj}|u_{J}|^2\ge 0\quad\T{ on $\partial D$},
\end{eqnarray}
 for any $u$  of degree $k\ge q>q_o$ (resp. $k\le q<q_o$) satisfying $u_{nK}|_{\partial D}=0$ for any $K$.
 \el
The proof is the same as in \cite{Z00}.

For convenient writing, we shall use the notation $A\lesssim B$ to mean $A\le cB $  for some constant $c$, which is independent of relevant parameters. And $A\cong B$ if $A\lesssim B$ and $B\lesssim A$.

\section{The basic estimates on q-pseudoconvexity/concavity}

In this section we prepare some inequalities which are needed for the subelliptic estimates of our Theorem~\ref{t1.1}. The key technical tool of our discussion are the so call Hormander-Kohn-Morrey estimates contained in the following proposition. 
Let $D$ be a domain with smooth boundary defined by $r=0$ in a neighborhood of $z_o$. Let $\om_1,...,\om_n=\partial r$ be an orhtonormal basis of $(1,0)$ forms and $L_1,...,L_n$ the dual basis of $(1,0)$ vector fields.

For         $0\le k\le n$, we write a general $k$-form $u$ as 
\begin{eqnarray*}
u={\sum_{|J|=k}}'u_J\bar{\om}_J,
\end{eqnarray*}
where $\sum'$ denotes summation restricted to ordered multiindices $J=\{j_1,...,j_k\}$ and where  $\bom_J=\bom_{j_1}\we...\we\bom_{j_k}$.  When the multiindex is no more ordered, it is understood that the coefficient $u_J$ is an antisymmetric function of $J$; in particular, if $J$ decomposes into $jK$, then $u_{jK}=\T{sign}\binom{J}{jK}u_J$. 
We  define $\langle{u,u}\rangle$ by 
$\langle{u,u}\rangle=|u|^2={\sum}'_{|J|=k}|u_J|$; this definition is independent of the choice of orthonormal basis $\om_1,...,\om_n$. The coefficients of our forms are taken in various spaces $\Lambda$ such as $C^\infty(\bar{D}), C^\infty({D}),C^\infty_c(\bar{D}), L^2(D)$ and the corresponding spaces of $k$-forms are denoted by $\Lambda^k$. Though our a priori estimates are proved over smooth forms, they are stated in Hilbert norms. Thus, let $||u||$ be the $H^0=L^2$ norm and, for a real function $\phi$, let the weighted $L^2$-norm be defined by 
$$\no{u}^2_{H^0_\phi}:={\sum_{|J|=k}}'\int_D e^{-\phi}|u_J|^2dv$$
where $dv$ is the element of volume in $\mathbb{C}^n$.
%%%%
We begin by noticing that $\bar\partial$ is closed, densely defined. Also,   its domain $D_{\bar\partial}$ certainly contains smooth forms and its action is expressed by    
 \begin{equation}
 \Label{3.1}
 \bar\partial u=\sumK\underset{\overset{ij=1,...,n}{i<j}}{\sum}(\bar L_iu_{jK}-\bar L_ju_{iK})\bar D _i\wedge\bar\omega_j\wedge\bar\omega_K+...,
 \end{equation}
 where dots denote terms in which no differentiation of $u$ occurs.  
 
Let $\bar\partial^*$ be the adjoint of $\bar\partial^*$. The operator $\bar\partial^*$ is still closed, densely defined but it is no more true that smooth forms belong to $D_{\bar\partial^*}$. For this, they must satisfy certain boundary conditions. Namely,  integration by parts shows that a form $u$ of degree $k$ cannot belong to $D_{\bar\partial^*}$  unless 
 $$
 \underset{j=1}{\overset n\sum}\int_{\partial D}e^{-\phi}L_j(r)u_{jK}\psi_K ds=0\quad\T{ for any $K$ and any $ \psi_K$ of degree $k-1$}.
 $$
This means that 
  $ \underset{j=1}{\overset n\sum} L_j(r)u_{jK}|_{\partial D}\equiv0\,\,\T{ for any } K$. (Here $ds$ is the element of hypersurface in $\partial D $.) Since we have chosen our basis with the property $L_j(r)|_{\partial D}=\kappa_{jn}$ (the Kronecker's symbol), we then conclude 
\begin{equation}
\Label{3.2}
 \T{ $u$ belongs to $ D_{\bar\partial^*}$ iff $u_J|_{\partial D}=0$ whenever $n\in J$.}
\end{equation}
We call {\it tangential } a form which belongs to $D_{\bar{\partial}^*}$.
Let  $\mathcal L^\phi_j$  be the formal $H^0_\phi$-adjoint of $-L_j$;
over a tangential form the action of  the Hilbert adjoint of $\bar\partial$, coincides with that of its ``formal adjoint" and is therefore expressed by a ``divergence operator":
 \begin{equation}
 \Label{3.3}
 \bar\partial^{*}_\phi u=-\sumK\underset j\sum\mathcal L_j^\phi (u_{jK})\bar\omega_K+...\quad 
 \T{ for any } u\in D_{\bar\partial^*},
 \end{equation}
 where dots denote an error term in which $u$ is not differentiated and $\phi$ does not occur.
By developing the equalities \eqref{3.1} and \eqref{3.3} by means of integration by parts, we get the proof of the following crucial result.

\bp
\Label{p3.1}
Let $D$ be a smoothly bounded domain  and fix arbitrarily an index $q_0$ with $0\le q_0\le n-1$. Then for a suitable  $C>0$ and any $u\in C^\infty(\bar{D})^k\cap D_{\bar{\partial}^*}$, we have 
\begin{eqnarray}
\Label{3.4}
\no{\bar{\partial} u}^2_{H^0_\phi}&+&\no{\bar{\partial}^*_\phi u}^2_{H^0_\phi}+C\no{u}^2_{H^0_\phi}\ge\\
&+&{\sum_{|K|=k-1}}'\sum^{n}_{i,j=1}\int_De^{-\phi}\phi_{ij}u_{iK}\bar{u}_{jK}dv-{\sum_{|J|=k}}'\sum^{q_o}_{j=1}\int_De^{-\phi}\phi_{jj}|u_{J}|^2dv\\
&+&{\sum_{|K|=k-1}}'\sum^{n-1}_{i,j=1}\int_{\partial D}e^{-\phi}r_{ij}u_{iK}\bar{u}_{jK}ds-{\sum_{|J|=q}}'\sum^{q_o}_{j=1}\int_{\partial D}e^{-\phi}r_{jj}|u_{J}|^2ds\\
&+&(1-\alpha)(\sum^{q_0}_{j=1}\no{\mathcal{L}^{\phi}_j u}^2_{H^0_\phi}+\sum^{n}_{j=q_o+1}\no{\bar{L_j}u}^2_{H^0_\phi}).
\end{eqnarray}
\ep
We refer for instance to \cite{Z08} for the proof of Proposition~\ref{p3.1}.
We note that there is no relation between $k$ and $q_o$ in above inequality and that $C$ and $\alpha$ are independent of $\phi$ (and $u$).  By choosing $\phi$ so that $e^{-\phi}$ is bounded, we may remove the weight functions in \eqref{3.4} to get some inequalities that are useful for the proof of Theorem~\ref{t1.1}. 
We write $Q(u,u)=\no{\bar{\partial} u}^2+\no{\bar{\partial}^* u}^2$.
\\
\bt
\Label{t3.1}
Assume that the hypotheses of Theorem~\ref{t1.1} be fulfilled. Then, for a suitable neighborhood $V$ of $z_o$ and for $\delta $ small, we have
\begin{eqnarray}
\Label{3.5}
\no{u}^2+\delta^{-2\epsilon}\int_{S_\delta}|u|^2dv +\sum^{q_o}_{j=1}\no{L_j u}^2+\sum^{n}_{j=q_o+1}\no{\bar{L}_ju}^2 \lesssim Q(u,u)
\end{eqnarray} 
for any $u\in C^{\infty}_c(\bar{D}\cap V)^k\cap D_{\bar{\partial}^*}$ with $k\ge q$ (resp. $k\le q$) when $q>q_o$ (resp. $q<q_o$).
\et

\bpf
We use twice Proposition~\ref{p3.1} and in both cases, owing to the assumption of $q$-pseudoconvexity (resp. $q$-pseudoconcavity)  we have the crucial fact that the boundary integrals are $\ge 0$ for any $k\ge q>q_o$ (resp. $k\le q<q_o$).
We first use Proposition~\ref{p3.1} under the choice $\phi\equiv 0$  and get 
\begin{eqnarray}
\Label{3.6}
Q(u,u)+C\no{u}^2\gtrsim \sum^{q_o}_{j=1}\no{ L_j u}^2+\sum^{n}_{j=q_o+1}\no{\bar{L}_ju}^2\quad  u\in C^{\infty}(\bar{D}\cap U)^k\cap D_{\bar{\partial}^*}. 
\end{eqnarray}
We use again Proposition~\ref{p3.1}, this time for $\phi=\chi(\phi^\delta)$.
In this case, the second line of \eqref{3.4} splits into two terms
\begin{multline}
\Label{3.8}
\int_De^{-\chi(\phi^\delta)}\dot\chi\Big(\sumK\underset{ij=1}{\overset n\sum}\phi^\delta_{ij}u_{iK}\bar u_{jK}-\underset{j=1}{\overset{q_o}\sum}\phi^\delta_{jj}|u|^2\Big)dv
\\
+\int_De^{-\chi(\phi^\delta)}\ddot\chi\Big(\sumK|\underset{j=1}{\overset n\sum}\phi^\delta_ju_{jK}|^2-\underset{j=1}{\overset{q_o}\sum}|\phi^\delta_j|^2|u|^2\Big)dv.
\end{multline}
We also have
\begin{eqnarray}
\Label{3.7}
\no{\bar{\partial}^*_{\chi(\phi^\delta)} u}^2_{H^0_{\chi(\phi^\delta)}}\le 2\no{\bar{\partial}^* u}^2_{H^0_{\chi(\phi^\delta)}}+2\sumK\no{\dot\chi^2\sum^n_j \phi^\delta_ju_{jK}}^2_{H^0_{\chi(\phi^\delta)}}.
\end{eqnarray}
Remark that $|\underset{j=1}{\overset{q_o}\sum}({\chi(\phi^\delta)})_j u|^2=|\dot\chi|^2|\underset{j=1}{\overset{q_o}\sum}\phi^\delta_j|^2|u|^2$. Thus we get from \eqref{3.4}, under the choice of the weight $\chi(\phi^\delta)$, and taking into account \eqref{3.8} and \eqref{3.7}:
\begin{equation}
\Label{3.9}
\begin{split}
\no{\bar{\partial} u}^2_{H^0_{\chi(\phi^\delta)}}&+\,2\no{\bar{\partial}^*_{\chi(\phi^\delta)}u}^2_{H^0_{\chi(\phi^\delta)}}+2C\no{u}^2_{H^0_{\chi(\delta)}}\\
&\ge\,\int_D\dot\chi e^{-\chi(\phi^\delta)}\Big({\sum_{|K|=k-1}}'\sum^{n}_{i,j=1}{\phi^\delta}_{ij}u_{iK}\bar{u}_{jK}dv-{\sum_{|J|=q_0}}'\sum^{q_o}_{j=1}{\phi^\delta}_{jj}|u_{J}|^2\Big)dv
\\
&+\, \int_D(\ddot\chi-2\dot\chi^2) e^{-\chi(\phi^\delta)}\sumK|\underset{j=1}{\overset n\sum}\phi^\delta_ju_{jK}|^2dv\,\,-\,\,\int_D\ddot\chi  e^{-\chi(\phi^\delta)}\underset{j=1}{\overset{q_o}\sum} |\phi^\delta_j|^2|u|^2 dv.
\end{split}
\end{equation}
We now specify our choice of $\chi$. First, we  want $\ddot\chi\ge 2\dot\chi^2$ so that the first sum in the third line can be disregarded. Keeping this condition, 
we need an opposite estimate which assures that the absolute value of the last negative term in the third line of 
\eqref{3.9} is controlled by one half of the second line. If $c$ is the constant of \eqref{1.5}, the above condition is fulfilled as soon as $\frac{2\ddot \chi}{\dot\chi}\leq  c$. If we then set $\chi:=\frac12e^{\frac c{2}(t-1)}$ then both requests are satisfied; (we also notice that $\dot \chi^2<<\dot \chi$ because $c<<1$). 
 Thus our inequality continues as
\begin{equation}
\begin{split}
\Label{3.10}
\quad\quad&\ge \,\frac{1}{2}\int_D\dot\chi e^{-\chi(\phi^\delta)}\Big({\sum_{|K|=k-1}}'\sum^{n}_{i,j=1}{\phi^\delta}_{ij}u_{iK}\bar{u}_{jK}dv-{\sum_{|J|=q_0}}'\sum^q_{j}{\phi^\delta}_{jj}|u_{J}|^2\Big)dv\\
&\geq\,\frac{1}{2}\int_{S_\delta} \dot\chi e^{-\chi(\phi^\delta)}\Big({\sum_{|K|=k-1}}'\sum^n_{i,j=1}{\phi^\delta}_{ij}u_{iK}\bar{u}_{jK}-{\sum_{|J|=k}}'\sum^q_{j}{\phi^\delta}_{jj}|u_{J}|^2\Big)dv\\
&\ge\,\delta^{-2\epsilon}\int_{S_\delta}\frac{c}{2}\dot\chi e^{-\chi(\phi^\delta)}|u|^2dv.
\end{split}
\end{equation} 
Here we are using the two main assumptions for our weights $\phi^\delta$, that is, 
\eqref{1.5}(with the right side replaced by $0$)  to get the second inequality and \eqref{1.6} as for the third.
Thus the first line of \eqref{3.8} is bigger or equal to the last of \eqref{3.10}.
We want to remove the weight from the resulting inequality. The first term can be handled owing to 
 $e^{-\chi(\phi^\delta)}\le 1$ on $\bar D\cap V$ and the second owing to $\dot \chi  e^{-\chi(\phi^\delta)}\ge c\ge0$ on $S_\delta\cap V$ which follows in turn from $|\phi^\delta|<1$. We end up with the unweighted estimate
\begin{eqnarray}
\Label{3.10bis}
\no{\bar{\partial} u}^2+\no{\bar{\partial}^* u}^2+C\no{u}^2
\simge\delta^{-2\epsilon}\int_{S_\delta}|u|^2dv.
\end{eqnarray}
Now, for fixed $\delta_o$ and for $V$ contained in the $\delta_o$-ball centered at $z_o=0$, the term $C||u||^2$ in the left of \eqref{3.10bis} can be absorbed in the right. Thus we end up with the estimate
\begin{eqnarray}
\Label{3.11}
\no{\bar{\partial} u}^2+\no{\bar{\partial}^* u}^2\simge\delta^{-2\epsilon}\int_{S_\delta}|u|^2dv+\no{u}^2
\end{eqnarray}
 for any $u\in C^{\infty}_c(\bar{D}\cap V)^k\cap D_{\bar\partial^*}$ and $\delta\le \delta_o$.\\

Combining \eqref{3.6} and \eqref{3.11}, we get \eqref{3.5} which concludes the proof of the theorem.

\epf

\section{Proof of Theorem~\ref{t1.1}}

Let $V$ be a neighborhood of a given point $z_o\in\partial D$, let $(t, r)$ be smooth coordinates in $V$  with $t=(t_1,...,t_{2n-1})$ and let $\tau$ be dual coordinates to $t$. For a function $u$ supported in $V$, one defines the tangential Fourier transform by 
$$\hat{u}(\tau, r)=\int_{\mathbb{R}^{2n-1}}e^{-it\tau}u(t,r)dt,$$
and  the tangential $H^s$-Sobolev norm  by 
$$|||u|||^2_s=\no{\Lambda^s u}^2 =\int^0_{-\infty}\int_{\mathbb{R}^{2n-1}}(1+|\tau|^2)^s|\hat{u}(\tau,r)|^2d\tau dr,$$
where $\Lambda^s$ is the tangential Bessel potential of order $s$. We note that when $s=0$ then $|||u|||_0=\no{u}$ is the usual $L^2$-norm.  We refer to \cite{FK72} for further details.

We remark that if $D_i$ is $\frac{\partial}{\partial t_j}$ or $\frac{\partial}{\partial r}$ then
 \begin{equation*}
 \begin{split}
||u||^2_s&=\sum_{i}|||D_iu|||^2_s
\\
&\cong |||u|||^2_{s+1}+|||D_ru|||^2_s.
\end{split}
\end{equation*}
The next result contains the key estimate in the proof of Theorem~\ref{t1.1}. 
\bl
\Label{l4.1}Let $U$ be a special boundary chart for $D$. Then for all $z_o\in \partial D\cap U$ there exsists a neighborhood $V\subseteq U$ of $z_o$ such that 
$$|||u|||_\epsilon^2\lesssim \sum_{j\le q_0}|||L_ju|||_{\epsilon-1}^2+\sum_{j\ge q_0+1}|||\bar{L}_j u|||_{\epsilon-1}^2+||u_b||^2_{\epsilon-\frac{1}{2}},~~~u\in C^\infty(V\cap \bar{D})^k\cap D_{\bar{\partial}^*}$$ 
where $u_b:=u|_{\partial D}$ and $\epsilon \le \frac{1}{2}$
\el

The above lemma is a variant of Theorem (2.4.5) of \cite{FK72}
to which we refer for the proof. Notice that
 on one hand our statement is more general because we choose any $\epsilon\leq\frac12$ instead of $\epsilon=\frac12$. On the other, we specialize the choice of a general elliptic system to the case of $\{L_j\}_{j\leq q_0}\cup\{\bar L_j\}_{q_0+1\leq j\leq n}$.\\

For the proof of Theorem~\ref{t1.1}, we use a method derived from \cite{C87}.
Let $p_k(t), k=0,1,...$ be a sequence of functions with $\sum^\infty_{k=0}p^2_k(t)=1$, $p_k(t)\equiv 0$ if $t\not\in(2^{k-1},2^{k+1})$ with $k\ge 1$ and $p_0(t)\equiv 0, t\ge 2$. We can also choose $p_k$ so that 
$$|p_k'(t)|\le C 2^{-k}.$$
Let $P_k$ denote the operator defined by 
$$(\widehat{P_ku})(\tau,r)=p_k(|\tau|)\hat{u}(\tau,r)$$
where $\hat{u}$ is the tangential Fourier transform. 
Let $\R^{2n}_-:=\{z:\,r(z)<0\}$ and denote by $\mathcal S(\R^{2n}_-)$ the Schwartz space of $C^\infty(\R^{2n}_-)$-functions which are rapidly decreasing at $\infty$.

\bl
\Label{l4.2}
For $f,u\in \mathcal S(\mathbb{R}_{-}^{2n})$ and $\sigma\in \mathbb{R}$ then 
$$\sum^\infty_{k=0}2^{2k\sigma}||[P_k, f]u||^2\lesssim |||u|||_{\sigma-1}^2.$$ 
\el

\bl
\Label{l4.3}
Let $L$ be a tangential vector field with coefficients in $C^\infty_0(\mathbb{R}^{2n}_-)$. Then 
$$\sum^{\infty}_{k=0} \no{[P_k,L]u}^2\le C\no{u}^2.$$
\el

The proof of Lemma 4.2 and Lemma 4.3 can be found in \cite{C87}. We remark that if we replace $u\in S(\mathbb{R}^{2n}_{-})$ by  $u\in C^\infty(\bar D\cap U)^k\cap D_{\bar\partial^*}$, then the two lemmas above still hold.\\

\noindent
{\it Proof of Theorem~\ref{t1.1}}
It suffices to prove the weaker version of \eqref{1.3,5} in which $||\cdot||_\epsilon$ is replaced 
by $|||\cdot|||_\epsilon$. In fact, $D_r$ can be expressed as a linear combination of $\bar L_n$ and a suitable ``totally real tangential" vector field that we denote by $ T $. 
We have
\begin{equation*}
\begin{cases}
Q(u,u)\simgeq ||\bar L_nu||^2,
\\
|||u|||^2_\epsilon\simgeq || T  u||^2_{\epsilon-1}.
\end{cases}
\end{equation*}
It follows 
\begin{equation*}
\begin{split}
||u||^2_\epsilon&=||D_r(u)||^2_{\epsilon-1}+|||u|||^2_\epsilon
\\
&\simleq Q(u,u)+|||u|||^2_\epsilon,
\end{split}
\end{equation*}
 which proves the claim. 
By Lemma 4.1 and Theorem~\ref{t3.1}, we get for any $u\in C^{\infty}(\bar{D}\cap V)^k\cap \Dom(\bar{\partial}^*)$ with $k\ge q+1$
\begin{eqnarray*}
|||u|||_\epsilon^2 &\lesssim& \sum_{j=0}^{q}|||L_j u|||^2_{\epsilon-1}+\sum_{j=q+1}^{n}|||\bar{L}_j u|||_{\epsilon-1}^2+\no{u_b}^2_{\epsilon-1/2}\\
&\lesssim& Q(u,u)+\no{u_b}^2_{\epsilon-1/2}.
\end{eqnarray*}
Now, we estimate $\no{u_b}^2_{\epsilon-1/2}$.
Let $\chi_k\in C^\infty_c(-2^{-k},0]$  with $0\leq\chi_k\leq1$ and $\chi_k(0)=1$. We have  the elementary inequality 
$$|g(0)|^2\le \frac{2^{k}}{\eta }\int^0_{-2^{-k}}|g(r)|^2dr+2^{-k}\eta\int^0_{-2^{-k}}|g'(r)|^2dr,$$
which holds for any $g$ such that $g(-2^{-k})=0$.
If we apply it for $g(r)=\chi_k(r)P_ku(\cdot,r)$, we get
\begin{eqnarray*}
\no{u_b}_{\epsilon-1/2}^2&\cong&\sum_{k=0}^{\infty}2^{2k(\epsilon-1/2)}\no{\chi_k(0)P_ku_b}^2\\
&\le&\eta^{-1} \sum_{k=0}^{\infty} 2^{2k\epsilon}\int_{-2^{-k}}^0\no{\chi_k P_ku(.,r)}^2dr+\eta\sum_{k=0}^{\infty}2^{2k(\epsilon-1)}\int_{-2^{-k}}^0 \no{D_r\Big(\chi_k P_ku(.,r)\Big)}^2dr\\
&=&\underbrace{\eta^{-1} \sum_{k=0}^{\infty} 2^{2k\epsilon}\int_{-2^{-k}}^0\no{\chi_k P_ku(.,r) }^2dr}_{I} +\underbrace{\eta\sum_{k=0}^{\infty}2^{2k(\epsilon-1)}\int_{-2^{-k}}^0 \no{D_r(\chi_k) P_ku(.,r)}^2dr}_{II}\\
&&+\underbrace{\eta\sum_{k=0}^{\infty}2^{2k(\epsilon-1)}\int_{-2^{-k}}^0 \no{\chi_k D_r\Big(P_ku(.,r)\Big)}^2dr}_{III}.
\end{eqnarray*}
Observe that $\chi_k\le 1$ and recall Theorem~\ref{t3.1}  that we apply for $P_ku$ and $\delta=2^{-k}$.  Thus the first sums above can be estimated by
\begin{eqnarray*}
(I) &\le& \eta^{-1} \sum_{k=0}^{\infty} 2^{2k\epsilon}\int_{-2^{-k}}^0\no{P_ku(.,r)}^2dr\\
 &\lesssim& \eta^{-1}\sum_{k=0}^{\infty} Q(P_ku,P_ku).
\end{eqnarray*}
We note that $Q(w,w)$ can be written as a finite sum of terms of the type   $$M_i=a_i T _i+b_iD_r+c_i,$$ where $ T _i$ are tangential vector fields. Hence
\begin{eqnarray*}
\sum_{k=0}^{\infty}Q(P_ku,P_ku)&\le& \sum_{k=0}^{\infty}\Big(\no{P_k\bar{\partial}u}^2+\no{P_k\bar{\partial}^*u}^2\Big)+\sum_{i}\sum_{k=0}^{\infty}\no{[M_i,P_k]u}^2\\
&\lesssim&Q(u,u)+\sum_{i}\sum_{k=0}^{\infty}\no{[a_iT_i,P_k]u}^2+\sum_{i}\sum_{k=0}^{\infty}\no{[b_i,P_k]D_r(u)}^2+|||u|||^2_{-1}
\\
&\lesssim& Q(u,u)+\no{u}^2+|||D_r(u)|||^2_{-1}, 
\end{eqnarray*}
where the estimates on the commutator terms follow by Lemma~\ref{l4.2} and Lemma~\ref{l4.3}. As it has already been remarked, $D_r(u)$ can be expressed as a linear combination of $\bar{L}_n u$ and $Tu$ for some tangential vector field $T$. Then
\begin{eqnarray*}
 |||D_r(u)|||^2_{-1}&\lesssim& |||\bar{L}_n u|||^2_{-1}+||| T u|||^2_{-1}\\
&\lesssim& \no{\bar{L_n}u}^2+||u||^2\\
&\lesssim& Q(u,u)
\end{eqnarray*}
where the last line follows from Theorem~\ref{t3.1}. 
 
We now estimate (II). Since $D_r(\chi_k)\le 2^k$, we get
$$(II)\le \eta\sum_{k=0}^{\infty}2^{2k\epsilon}\int_{-2^{-k}}^0 \no{P_ku(.,r)}^2dr\le \eta\sum_{k=0}^{\infty}2^{2k\epsilon} \no{P_ku}^2\cong \eta|||u|||_\epsilon^2.$$
As for the term (III), we have $D_rP_k=P_kD_r$ and $\chi_k\le 1$. Also $D_r=a\bar{L}_n+b T $ as before. Thus 
\begin{eqnarray*}
 (III)&\le& \eta\sum_{k=0}^{\infty}2^{2k(\epsilon-1)}\no{P_kD_r(u)}\cong \eta|||D_r(u)|||_{\epsilon-1}\\ &\lesssim& \eta\left(|||\bar{L}_n u|||^2_{\epsilon-1}+||| T u|||^2_{\epsilon-1}\right)\\
&\lesssim&\eta Q(u,u)+\eta|||u|||_{\epsilon}^2.
\end{eqnarray*}
  Combining all our estimates of $\no{u_b}_{\epsilon-1/2}$, we obtain
$$\no{u_b}_{\epsilon-1/2}\lesssim \eta^{-1}Q(u,u)+\eta |||u|||_{\epsilon}.$$
Summarizing up, we have shown that
$$|||u|||_\epsilon\lesssim \eta^{-1}Q(u,u)+\eta|||u|||_{\epsilon}^2.$$
Choosing $\eta>0$ sufficiently small, we can move the term $\eta|||u|||^2_\epsilon$ into the left-hand-side and get $$|||u|||_\epsilon^2\lesssim Q(u,u).$$ The proof is complete.

\hskip12cm $\Box$

\section{Proof of  Theorem~\ref{t1.2}}

We note that $r=2x_n+h$ is a graphing function and denote by $z\mapsto z^*$ the projection $\bar D\to \partial D$ along the $x_n$-axis. We denote by $\partial r^\perp(z), \,z\in V$ the bundle orthogonal to $\partial r=\omega_n$ and note that $\partial r^\perp|_{\partial D}=T^{1,0}\partial D$. We have the evident equalities
\begin{equation}
\Label{5.1}
\begin{cases}
(r_{ij}(z))_{ij=1}^{n-1}=(r_{ij}(z^*))_{ij=1}^{n-1},
\\
\partial r^\perp(z)=\partial r^\perp(z^*).
\end{cases}
\end{equation}
Thus \eqref{5.1} relates $L_r|_{T^\C \partial D}$ on $\partial D\cap V$  to $L_r|_{\partial  r^\perp}$ on the whole of $\bar D\cap V$; in particular, \eqref{1.3} passes from $\partial D\cap V$ to the whole of $\bar D\cap V$. 
Since $\tilde r$ is obtained by adding to $r$ terms which in turn satisfy \eqref{1.3} in the two respective cases, then one can prove that $\tilde D$ is $q$-pseudoconvex (resp. $q$-pseudoconcave) in the sense of its ``exhaustion" functions (though this is not clear for defining functions). We do not enter into these details and just show, in the beginning of the proof that \eqref{5.1}, which concerns the behavior of $r$ on $\partial D$, turns into a similar property of a weight $\phi$ in $\bar D$.

We choose a local basis  $\omega_1,...,\omega_n=\partial \tilde r$  of $(1,0)$-forms and denote by $L_1,...,L_n$ the dual basis of $(1,0)$-vector fields; we may assume that $L_j(z_o)=\partial_{z_j}$. 
Thus, by an orthonormal change in the system $\T{Span}\{L_1,...,L_{n-1}\}$, we can assume that \eqref{1.3} is satisfied on $\partial D$. 
We now construct the weight $\phi$ which satisfies the assumptions of Theorem~\ref{t1.1}; for this we distinguish $q>q_o$ from $q<q_o$.
\\
\underline{The case $q$-pseudoconvex.}
We set for  a suitable constant $\lambda>0$
\begin{equation}
\Label{5.2}
\psi=-\log(-\tilde r+\delta)+\lambda|z|^2+\underset{j=q}{\overset{n-1}\sum}\log(|z_j|^2+\delta^{\frac1{m_j}}),
\end{equation}
and define $\phi:=c|\log\,\delta|^{-1}\psi$ where $c$ is an irrelevant constant needed to get the bound $1$ in \eqref{1.5} and \eqref{1.6}.
We set $\psi^I=-\log(-\tilde r+\delta)+\lambda|z|^2$ and denote by $\psi^{II}$ the remaining term in the right of \eqref{5.2}; thus $\psi=\psi^I+\psi^{II}$. We have
\begin{equation}
\Label{5.3}
\begin{split}
\psi_{ij}^I&=(-\tilde r+\delta)^{-1} \tilde r_{ij}+\lambda\kappa_{ij}
\\
&=(-\tilde r+\delta)^{-1}r_{ij}\,+\,\lambda\kappa_{ij}+(-\tilde r+\delta)^{-1}(\partial_{z_j}\partial_{\bar z_j}h_j)\kappa_{ij}\,+\,\mathcal E\quad\T{for $i,j\le n-1$},
\end{split}
\end{equation}
where $\mathcal E$ is an error of type $\mathcal E=O(|z|)(-\tilde r +\delta)^{-1}\sum_j(\partial_{z_j}\partial_{\bar z_j}h_j)$. 
We also have
\begin{equation}
\Label{5.3bis}  
\psi_{nn}^I=(-\tilde r+\delta)^{-2}.
\end{equation}
(where $\kappa_{ij}$ continues to denote the Kronecker's symbol) and
\begin{equation}
\Label{5.4}
\psi_{ij}^{II}=\Big(\frac{\delta^{\frac1{m_j}}}{(|z_j|^2+\delta^{\frac1{m_j}})^2}\Big)\kappa_{ij}
\end{equation}
When taking $\underset{ij}\sum\cdot-\underset{j=1}{\overset{q_o}\sum}\cdot$
of $(-\tilde r+\delta)^{-1}r_{ij}\,+\,\lambda\kappa_{ij}$ from \eqref{5.3} and of $(-r+\delta)^{-2}$ from \eqref{5.3bis} the result is $\ge0$. This is true for $(r_{ij}(z))_{ij}|_{\partial r^\perp(z)}$ on account of \eqref{5.1}. But what is left is just
$$
\lambda|\omega|^2+(-\tilde r+\delta)^{-2}|\partial r|^2+(-\tilde r+\delta)^{-1}2\Re\underset{j=1}{\overset{n}\sum}r_{nj}\partial r\otimes\bar\omega_j,
$$
which is positive.
We also discard all terms of type $(\partial_{z_j}\partial_{\bar z_j}h_j)\kappa_{ij}$ and $\delta^{\frac1{m_j}}\kappa_{ij}$ for $i$ or $j$ $\le k-1$ 
in addition to $\mathcal E$ because they can be made positive by adding a small amount of terms for which $i,j\ge k$ on account of the estimates \eqref{5.6} and \eqref{5.7} which follow.
For the remaining terms $(\partial_{z_j}\partial_{\bar z_j}h_j)$, we note that we have $(\partial_{z_j}\partial_{\bar z_j}h_j)\simge|z_j|^{2m_j-2}$. We end up with the estimate
\begin{multline}
\Label{5.5}
\sumK\underset{ij=1}{\overset n\sum}\psi_{ij}u_{iK}\bar u_{jK}-\underset{j=1}{\overset {q_o}\sum}\psi_{jj}|u|^2
\\
\geq \underset{j=k}{\overset {n-1}\sum}\Big((-\tilde r+\delta)^{-1}|z_j|^{2m_j-2}+\frac{\delta^{\frac1{m_j}}}{(|z_j|^2+\delta^{\frac1{m_j}})^2}\Big)\sumK|u_{jK}|^2+(-r+\delta)^{-2}\sumK|u_{nK}|^2.
\end{multline}
We now inspect the coefficients in the right of \eqref{5.5}.
First, let $z\in S_\delta$, that is, $-r>\delta$. Given a coefficient $u_J$ of $u$, the index $J$ contains for sure at least one $j$ such that $k\le j\le n-1$ and thus 
$u_J=\T{sign}\binom J{jK} u_{jK}$ for a suitable $K$. If, for this $j$, $|z_j|^2\ge \delta^{\frac1{m_j}}$, then
\begin{equation}
\Label{5.6}
(-\tilde r+\delta)^{-1}|z_j|^{2m_j-2}\simge\delta^{-\frac1{m_j}}.
\end{equation}
On the contrary, if $|z_j|^2\leq \delta^{\frac1{m_j}}$, then 
\begin{equation}
\Label{5.7}
\frac{\delta^{\frac1{m_j}}}{(|z_j|^2+\delta^{\frac1{m_j}})^2}\simge \delta^{-\frac1{m_j}}.
\end{equation}
In both cases, the terms in the left are $\ge \delta^{-2\epsilon_k}$ since $-\frac1{m_j}\le-\frac1{m_k}=-2\epsilon_k$. By combining \eqref{5.6} with \eqref{5.7}, we get the second of \eqref{1.6} for $\epsilon=\epsilon_k$. On the other hand, for any $j\le q_o$, we have $\frac{r_j}{-r+\delta}=0$ 
and therefore $\psi_j=\lambda O(|z|)$ which is estimated by $\bar\partial\partial(\lambda|z|^2)$. On the other hand,
 $\underset{ij}\sum\cdot-\underset{j=1}{\overset {q_o}\sum}\cdot$ is always $\ge0$ all over $\bar D\cap V$. This proves the second inequality in \eqref{1.5}. 

Finally, a normalization by a factor $c|\log\delta|^{-1}$ makes the weight bounded as required by the first of \eqref{1.5} and \eqref{1.6}, at the expenses of passing from $\delta^{-2\epsilon_k}$ to $\frac{\delta^{-2\epsilon_k}}{|\log\delta|}$ in \eqref{1.6}. 
Thus the weight $\psi$ satisfies all the requirements of Theorem~\ref{1.1} for any $\epsilon<\epsilon_k$ which implies subelliptic estimates of the corresponding order.
Incidentally, we notice that when $\epsilon_k=\frac12$, the term $\psi^{II}$ is needless and we can take a different normalization by defining $\phi=-\log\left(\frac{-r+\delta}{2\delta}\right)$; thus we get an even $\delta^{-1}$ on the right of \eqref{1.6}. 
For $\epsilon_k=\frac12$, a similar argument applies also to  the case $q$-pseudoconcave which follows and we will not insist on it.
\\
\underline{The case $q$-pseudoconcave.}
 We now define
 \begin{equation*}
 \psi=-\log(-\tilde r+\delta)\,-\,\lambda|z|^2\,+\,\underset{j=1}{\overset{k+1}\sum}
 \log(-\log(|z_j|^2+\delta^{\frac1{m_j}}))
 \end{equation*}
 where we point out the attention to the double $\log$.
Comparing with the case $q$-pseudoconvex, there is now an extra difficulty  for the weight to satisfy \eqref{1.5} (whereas \eqref{1.6} remains substantially unchanged)
because we do not have any longer $\phi_j=0$ for $j\le q_o$.
We write $\psi=\psi^I+\psi^{II}$ in the same way as in the previous case and will eventually define $\phi$ by a normalization $\phi=c\left|\log\delta\right|^{-1}\psi$. We have the analogous of \eqref{5.3} and \eqref{5.4} with the suitable sign. We apply $\underset{ij}\sum\cdot-\underset{j=1}{\overset{q_o}\sum}\cdot$ to $\psi^I+\psi^{II}$. 
When taking $\sum_{ij}\cdot-\underset{j=1}{\overset {q_o}\sum}\cdot$
we discard the contribution of 
$(-\tilde r+\delta)^{-1}r_{ij}\,+\,\lambda\kappa_{ij}$ in addition to the normal term $(-\tilde r+\delta)^{-2}$ because this contribution is positive as before. 
 We discard  the error term $\mathcal E$ because it can be made positive by the aid of a small amount of the remainder. This argument is the same as for the case $q$-pseudoconvex. 
 What we are left with is
\begin{multline}
\Label{5.8}
\sum_{ij}\cdot-\underset{j=1}{\overset {q_o}\sum}\cdot
\geq \underset{j=1}{\overset{k+1}\sum}\left((-r+\delta)^{-1}|z_j|^{2m_j-2}+\frac{-\delta^{\frac1{m_j}}}{(|z_j|^2+\delta^{\frac1{m_j})^2}}\frac1{|\log(|z_j|^2+\delta^{\frac1{m_j}})|}\right.
\\
+\,\left.\frac{|z_j|^2}{(|z_j|^2+\delta^{\frac1{m_j})^2}}\frac1{|\log(|z_j|^2+\delta^{\frac1{m_j}})|^2}\right)(|u|^2-\sumK|u_{jK}|^2).
\end{multline}
We write the coefficient in the right of \eqref{5.8} as $(A_j+B_j+C_j)$. The two first terms serve to get \eqref{1.6}, the third for \eqref{1.5}. (This latter was discarded as $\ge0$ in the case $q$-pseudoconvex; here it is essential because $\phi_j\neq0$ for $j\le q_o$). Reasoning as in the first half of the proof we get, for any $j\le k+1$
\begin{equation}
\Label{5.9}
A_j+B_j\simge \delta^{-\frac1{m_j}}\ge \delta^{-2\epsilon_k}\quad\T{ on $S_\delta\cap V$},
\end{equation}
 because $-\frac1{m_j}\le-\frac1{m_{k+1}}=-2\epsilon_k$ for any $j\le k+1$, along with
 \begin{equation}
 \Label{5.10}
 A_j+B_j\ge0\quad\T{ on $\bar D\cap V$.}
 \end{equation}
 We make the crucial remark for the case of concavity. If the degree of $u$ is $k$, then
 \begin{equation}
 \Label{5.11}
 \underset{j=1}{\overset{k+1}\sum}\left(|u|^2-\sumK|u_{jK}|^2\right)\ge |u|^2.
 \end{equation}
 From \eqref{5.9} and \eqref{5.11} we get \eqref{1.6}.
  We now need to prove that on $\bar D\cap V$ and for a suitable $\epsilon$ we have
  \begin{equation}
  \Label{5.12}
  \begin{split}
\underset{j=1}{\overset{k+1}\sum}\frac1{\log^2\cdot}\frac{|z_j|^2}{(|z_j|^2+\delta^{\frac1{m_j}})^2}&\left(|u|^2-\sumK|u_{jK}|^2\right)
\\
  &\ge \epsilon\underset{j=1}{\overset{k+1}\sum}\frac{|z_j|^2}{\log^2\cdot\,(|z_j|^2+\delta^{\frac1{m_j}})^2}|u|^2
  \\
&=\epsilon\underset{j=1}{\overset{k+1}\sum}|\phi_j|^2|u|^2.
  \end{split}
  \end{equation}
This would conclude the proof of \eqref{1.5}.
The last sum $\underset{j=1}{\overset{k+1}\sum}\cdot$ can be replaced by $\underset{j=1}{\overset{q_o}\sum}\cdot$ since $\phi_j=0$ for $j=k+2,...,q_o$. 
 Also, remember here that $\psi^I_j=0$ for any $j$ and $\psi^{II}_j=0$ for any $j\ge k+2$; this justifies the last equality in \eqref{5.12} which is true.
However, the first inequality is wrong. To make it true, we need a small perturbation of $\psi$. We take a vector $v$ 
in the unit sphere $S^k$ 
outside the first quadrant, set $\psi^{II\,v}:=\underset{j=1}{\overset{k+1}\sum}\log(-\log(|z_j|^2+\delta^{\frac1{m_j}})^{v_j})$, leave $\psi^I$ unchanged and define a new $\psi$ by
$$
\psi:=\psi^I+\frac12(\psi^{II}+\psi^{II\,v}).
$$
Inequalities \eqref{5.9} and \eqref{5.10} are stable under perturbation and thus will remain true for this new $\psi$. As for the first of \eqref{5.12}, we consider the vector field
$$
w(z):=\left(\frac{|z_j|}{\log^2\cdot\,(|z_j|^2+\delta^{\frac1{m_j}})}\right)_{j=1,...,k+1}.
$$
We also define
$$
\mu(z)= \frac{w(z)}{|w(z)|},\qquad \nu(z)=(\nu_j v_j)_{j=1,...,k+1};
$$
thus $|\mu|=1$ and $|\nu|\le 1$. Finally, we set
$$
u=\frac{(|u_{jK})_{j=1,...,k+1}}{\sum_j|u_{jK}|^2}.
$$
It suffices to prove that
$$
\frac12\left(\langle \mu,u\rangle^2+\langle \nu,u\rangle^2\right)\le 1-\epsilon.
$$
Now, we begin by noticing that
\begin{equation}
\Label{5.15}
\begin{cases}
\langle\mu,u\rangle\le 1,
\\
\langle\nu,u\rangle\le 1,
\end{cases}
\end{equation}
by Cauchy-Scwhartz inequality. Also, if the first of \eqref{5.15}  happens to be equality, that is, $\mu$ is parallel to $u$, then
\begin{equation*}
\begin{split}
\langle \nu,u\rangle&=\sum_jv_j\mu_ju_j
\\
&=\sum_jv_j\mu_j^2.
\end{split}
\end{equation*}
But for this to be $1$ we need both $\sum_j\mu_j^4=1$ and $v$ parallel to $(\mu^2_j)_{j=1,...,k+1}$ . If the first occurs then, since $\sum_j\mu_j^2=1$, we  have  $(\mu^2_j)=(\mu_j)$ (and both coincide with a cartesian vector): thus $(\mu^2_j)$ is not parallel to $v$. 
In conclusion if the first of \eqref{5.15} is equality, the second is not.
Therefore, the function 
$(u,\mu)\mapsto \frac12(\langle u,\mu\rangle^2+\langle u,\nu\rangle^2)$ has a minimum $<1$, say $1-\epsilon$, for $u\,\in S^k$ (and for $\nu=(\mu_jv_j)$).

\hskip14cm$\Box$

\section{Proof of Theorem~\ref{t1.3}}
\Label{s6}
\bl
\Label{l6.1}
We have
\begin{equation}
\Label{6.1}
\int_0^\delta...\int_0^\delta \frac{dx_1dy_1...dx_pdy_p}{(t\sum_{j=1}^{p}|t^{-\varepsilon}z_j|^{2m_j}+1)^{s}}\cong t^{-\sum_{1}^{p}\frac{1}{m_j}+2p\varepsilon},
\end{equation}
provided that $s>\frac1{m_1}+...+\frac1{m_p}+1$.
\el
\bpf
We can assume that $m_1\le m_2\le...\le m_p$. Put $a(t)=t\sum_{j=2}^{p}|t^{-\varepsilon}z_j|^{2m_j}+1$. First, we perform integration
$$M(z_2,...,z_p)=\int_0^\delta\int_0^\delta \frac{dx_1dy_1}{(t|t^{-\varepsilon}z_1|^{2m_1}+a(t))^{s}}.
$$ We also make a change of variables $z_1'=t^{\frac{1}{2m_1}-\varepsilon}a(t)^{-\frac{1}{2m_1}}z_1$ and get
$$M(z_2,...,z_p)=a(t)^{-s+\frac{1}{m_1}}t^{-\frac{1}{m_1}+2\varepsilon}\int_0^{t^{\frac{1}{2m_1}-\varepsilon}a(t)^{-\frac{1}{2m_1}}\delta}\int_0^{t^{\frac{1}{2m_1}-\varepsilon}a(t)^{-\frac{1}{2m_1}}\delta} \frac{dx'_1dy'_1}{(|z_1'|^{2m_1}+1)^{s}}.
$$
Since $$t^{\frac{1}{2m_1}-\varepsilon}a(t)^{-\frac{1}{2m_1}}\delta=\Big(\frac{t^{1-2\varepsilon m_1}\delta^{2m_1}}{\sum_{j=2}^p t^{1-2\varepsilon m_j}|z_j|^{2m_j}+1}\Big)^{\frac{1}{2m_1}}\ge C>0,$$
then
$$M(z_2,...,z_p)\cong a(t)^{-q+\frac{1}{m_1}}t^{-\frac{1}{m_1}+2\varepsilon}.$$
In conclusion, the left hand side of \eqref{6.1} is equivalent to 
$$t^{-\frac{1}{m_1}+2\varepsilon}\int_0^\delta...\int_0^\delta \frac{dx_2dy_2...dx_pdy_p}{(t\sum_{j=2}^{p}|t^{-\varepsilon}z_j|^{2m_j}+1)^{q-\frac{1}{m_1}}}.$$
Repetition of this argument for $z_2,...z_p$ yields the proof of the lemma. 

\epf

\noindent
{\it {Proof of Theorem~\ref{t1.3} (i)}} 
Let 
\begin{equation*}
\begin{cases}
\omega_j=dz_j-r_{\bar z_j}dz_n,
\\
\omega_n=\partial r
\end{cases}
\end{equation*}
be a  basis of $(1,0)$ forms. We note that for a $k$-form $u$ we have $u\in D_{\bar\partial*}$ if and only if its coefficients satisfy $u_{nK}|_{\partial D}\equiv0$ for any $|K|=k-1$. 
Let $L_j$ be the dual basis of $(1,0)$ vector fields; these are a perturbation of $\partial_{z_j}-r_{z_j}\partial_{z_n}\,\,j=1,...,n-1$ and $\underset{j=1}{\overset n\sum}r_{z_j}\partial_{z_j}$. 
We have 
\begin{itemize}
\item
$(\omega_{iK},\omega_{jK})=\kappa_{ij}+r_{z_i}r_{\bar z_j}\quad \text{for any }i,j\le n-1$,
\item
$(\omega_{jK},\omega_{nK})=0\quad \text{ for any }j\le n-1$,
\item
$(\omega_I,\omega_J)=0\quad \T{if $|I\cap J|\le k-2$}$,
\item
$
\left(\bar\partial^*u\right)_K=
\underset{j=1}{\overset n\sum}\underset{\{J:\,|J\cap jK|=k\}}\sum L_j(u_J)+\underset{j=1}{\overset{n-1}\sum}\underset {\{J:\,|J\cap jK|=k-1\}}\sum L_j(u_J)\left(O(r_{z_j})+\underset{i\in J}\sum O(r_{z_i})\right)$
\newline
\hskip3cm${}\quad\qquad+\,\T{error}$,
\end{itemize}
where ``error" denotes a term where no derivatives of $u$ occur.
We will deal with the form
$$
u_t=U_t\bar\omega_1\wedge...\wedge\bar\omega_k,
$$
where $U_t$ is  a functions which will be specified later.
We have for this form
\begin{equation*}
\begin{cases}
\bar\partial u\simeq \underset{j=k+1}{\overset n\sum}\bar L_j(U_t)\bar\omega_j\wedge\bar\omega_1\wedge...\wedge\bar\omega_k\,+\T{error},
\\
\begin{split}
\bar\partial^*u=&\underset{j=1}{\overset k\sum}L_j(U_t)\bar\omega_1\wedge...\overset\wedge{\bar\omega_j}\wedge...\wedge\bar\omega_k + \underset{j=1}{\overset{n-1}\sum}\underset{H\cap\{1,...,k,j\}\neq\emptyset}\sum L_j(U_t)\left(O(r_{z_j})
+\underset{i\le k}\sum O(r_{z_i})\right)\bar\omega_H
\\
&+\,\T{error}.
\end{split}
\end{cases}
\end{equation*}
In particular
\begin{equation}
\Label{6.3}
\begin{split}
||\bar\partial u||^2+||\bar\partial^*u||^2&\simleq \underset{j=k+1}{\overset n\sum}||\bar L_j U_t||^2+\underset{j=1}{\overset k\sum}||L_jU_t||^2
\\
&
+{\underset{i=k+1,...,n-1}\sum}||\left(O(|r_{z_i}|)+\underset{i\le k}\sum O(|r_{z_j}|)\right)L_iU_t||^2+||U_t||^2.
\end{split}
\end{equation}
We now set $U_t=f_t(z_n)\Phi_t(z)$ where 
\begin{equation*}
\begin{cases}
f_t(z',z_n)=(z_n-Q(z')-\frac{1}{t})^{-p}\\
\Phi_t(z)=\big(\Pi_{j=1}^{n-1}\phi(t^{\epsilon_k}x_j)\phi(t^{\epsilon_k}y_j)\big)\lambda(x_n)\phi(y_n).
\end{cases}
\end{equation*}
Here $\phi\in C^\infty_0(\mathbb{R})$ satisfies
\begin{equation*}
\phi(x)=\begin{cases}1 & x\le \delta\\
0 & x\ge 2\delta,
\end{cases}
\end{equation*}
where
$\delta$ is a small parameter, and $\lambda\in C^\infty_0(\mathbb{R})$ will be chosen later.\\
Since
\begin{equation*}
\begin{cases}
L_j(f_t)=0\quad\T{for any $j\le q_o$},
\\
\bar L_j(f_t)=0\quad\T{for $q_o+1\le j\le n-1$},
\\
\partial_{z_j}(f_t)=\partial_{\bar z_j}(f_t)=0\quad\T{for $j=q_o+1,...,n-1$},
\end{cases}
\end{equation*}
we can restrict the first sum in \eqref{6.3} to $j=n$ and the second to $j=q_o+1,...,k$; thus we get
\begin{equation}
\Label{6.4}
\begin{split}
Q(u_t,u_t)&\simleq \underset{ij=1}{\overset{q_o}\sum}||r_{z_i}\partial_{\bar z_j}(f_t)\Phi_t||^2+\underset{j=q_o+1}{\overset k\sum}||r_{z_j}\partial_{z_n}(f_t)\Phi_t||^2+\underset{j=1,...,k,i}{\underset{i=k+1,...,n-1}\sum}|||r_{z_j}||r_{z_i}|\partial_{z_n}(f_t)\Phi_t||^2
\\
&+\underset{ij=1}{\overset {n-1}\sum}||O^2(|r_{z_i}|)\partial_{z_j}f_t\Phi_t||^2+\underset{j=1}{\overset n\sum}||f_t\partial_{z_j}\Phi_t||^2+||U_t||^2.
\end{split}
\end{equation}
To estimate the first three sums in \eqref{6.4} we need to evaluate $r_{z_i}$ for $i=1,...,q_o$, next $r_{z_j}$ for $j=q_o+1,...,k$ and finally $r_{\bar z_i}r_{z_j}r_{z_i}$ for $i=k+1,...,n-1$, $j=1,...,k$. We perform the change of variables
\begin{equation*}
\begin{cases}
\tilde z_j=t^{\epsilon_k}z_j,\quad j\le n-1,
\\
\tilde z_n=tz_n.
\end{cases}
\end{equation*}
For $|\tilde z|\le 1$ we have for the first terms 
\begin{equation}
\Label{6.5}
\begin{split}
|r_{z_j}(z)|&=|z_j|^{4m_j-2}
\\
&=t^{\epsilon_k(4m_j-2)}\le t^{-2+2\epsilon_k}, \,\,j=q_o+1,...,k,
\end{split}
\end{equation}
where the last inequality follows from $m_j\ge m_k$. For the second terms we have
\begin{equation}
\Label{6.6}
\begin{split}
|r_{z_i}(z')|^2&\le|z'|^{4m-2}
\\
&=t^{\epsilon_k(4m-2)}\le t^{-2+2\epsilon_k},\,\,i=1,...,q_o,
\end{split}
\end{equation}
where the last estimate follows from $m\ge m_{q_o+1}\ge m_k$. For the third terms we extend the definition of $m_j$ to  $j\le q_o$ by putting $m_j=m$. We have, for $i\ge k+1,\,j\le k$ or $j=i$
\begin{equation}
\Label{6.7}
\begin{split}
|r_{z_j}(z)|^2|r_{z_i}(z)|^2&\le t^{-\epsilon_k(4m_i+4m_j-4)}
\\
&\le t^{-2-\epsilon_k(4m_j-4)}\le t^{-2},
\end{split}
\end{equation}
where the second inequality follows from $m_i\le m_k$. If we pass to estimate the terms in the second sum of \eqref{6.4}, we then have 
\begin{eqnarray*}
\no{r_{z_j}\frac{\partial f_t}{\partial z_n}\Phi_t}^2&\cong& \int |r_{z_j}|^2|z_n-Q(z')-\frac{1}{t}|^{-2p-2}\Phi^2_t(z)dx_1dy_1....dx_ndy_n \\
&\lesssim& \int \frac{|z_j|^{4m_j-2} }{\Big((\frac{1}{t}+Q(z')-x_n)^2+y_n^2\Big)^{p+1}}\Phi^2_t(z)dx_1dy_1....dx_ndy_n
\\
&\lesssim& t^{2p-2+2\epsilon_k-2(n-1)\epsilon_k}I_t,
\end{eqnarray*}
where 
$$
I_t=\int \frac{\big(\Pi_{j=1}^{n-1}\phi(\tilde{x}_j)\phi(\tilde{y}_j)\big)^2\lambda(t^{-1}\tilde{x}_n)^2\phi(t^{-1}\tilde{y}_n)^2}{\Big((1+tQ(t^{-\epsilon_k}\tilde{z}')-\tilde{x_n})^2+\tilde{y}_n^2\Big)^{p+1}}d\tilde{x}_1d\tilde{y}_1...d\tilde{x_n} d\tilde{y}_n.
$$
\\
We now perform integration in  $\tilde{y}_n$ from $-\infty$ to $+\infty$ and get
$$
I_t\lesssim \int \frac{\big(\Pi_{j=1}^{n-1}\phi(\tilde{x}_j)\phi(\tilde{y}_j)\big)^2\lambda(t^{-1}\tilde{x}_n)^2}{(1+tQ(t^{-\epsilon_k}\tilde{z}')-\tilde{x_n})^{2p+1}}d\tilde{x}_1d\tilde{y}_1...d\tilde{x}_{n-1}d\tilde{y}_{n-1}d\tilde{x_n}.
$$
\\
Next, we  integrate in $\tilde{x_n}$ from $$-\infty \text{ to } \Big(tQ(t^{-\epsilon_k}\tilde{z}')-t\sum_{j=q_0+1}^{n-1}|h_j(t^{-\epsilon_k}\tilde{z}_j)|^{2}\Big)/2,
$$ 
and get
\begin{eqnarray*}
I_t&\lesssim& \int \frac{\big(\Pi_{j=1}^{n-1}\phi(\tilde{x}_j)\phi(\tilde{y}_j)\big)^2}{(tQ(t^{-\epsilon_k}\tilde{z}')+t\sum_{j=q_0+1}^{n-1}|h_j(t^{-\epsilon_k}\tilde{z}_j)|^{2}+2)^{2p}}d\tilde{x}_1d\tilde{y}_1...d\tilde{x}_{n-1}d\tilde{y}_{n-1}\\
&\lesssim&\int_0^{2\delta}...\int_0^{2\delta}\frac{d\tilde{x}_{k+1}d\tilde{y}_{k+1}...d\tilde{x}_{n-1}d\tilde{y}_{n-1}}{(t\sum_{j=k+1}^{n-1}|t^{-\epsilon_k}z_j|^{2m_j}+1)^{2p}}\\
&\lesssim&t^{-\sum_{j=k+1}^{n-1}\frac{1}{m_j}+2(n-k-1)\epsilon_k}
\end{eqnarray*}
where the last inequality follow by Lemma~\ref{l6.1}.\\
In conclusion we have obtained
\begin{equation}
\Label{6.8}
\no{r_{z_j}\frac{\partial f_t}{\partial z_n}\Phi_t}^2\lesssim t^{2p-2+2\epsilon_k-2k\epsilon_k-\sum_{j=k+1}^{n-1}\frac{1}{m_j}}.
\end{equation}
The same integration combined with \eqref{6.7} yields the same estimate as \eqref{6.8} also for the terms $|||r_{z_j}||r_{z_i}|\partial_{z_n}(f_t)\Phi_t||^2$ for $i\ge k+1$ and $j\le k$. As for the terms in the first sum in \eqref{6.4}
with $i,j=1,...,q_0$, we have  
\begin{eqnarray*}
\no{r_{z_j}\frac{\partial f_t}{\partial \bar{z}_j}}^2&\cong&\int|r_{z_j}|^4|z_n-Q(z')-\frac{1}{t}|^{-2p-2}\Phi^2_t(z)dx_1dy_1....dx_ndy_n \\
&\lesssim& \int \frac{|z'|^{4m_k-2} }{\Big((\frac{1}{t}+Q(z')-x_n)^2+y_n^2\Big)^{p+1}}\Phi^2_t(z)dx_1dy_1....dx_ndy_n\\
&\lesssim& t^{2p-2+2\epsilon_k-2k\epsilon_k-\sum_{j=k+1}^{n-1}\frac{1}{m_j}}
\end{eqnarray*}
where the last inequality follows by the same technique as above.
 \\
By the same argument all the sums $\underset{ij=1}{\overset {n-1}\sum}||O^2(|r_{z_i}|)\partial_{z_j}f_t\Phi_t||^2$, the terms $\no{f_t\frac{\partial \Phi_t}{\partial z_j}}^2\,\,j=1,...,n$ and $\no{U_t}^2$ 
have the same estimate in terms of $t$. Combining all these estimates, we get the basic estimate from above for $Q(u_t,u_t)$
\begin{eqnarray}
\Label{6.9}
Q(u_t,u_t)&\lesssim& t^{2p-2+2\epsilon_k-2k\epsilon_k-\sum_{j=k+1}^{n-1}\frac{1}{m_j}}.
\end{eqnarray}
To calculate $|||u|||_{\epsilon}$ we use the boundary coordinates $(x_1,...,x_{n-1}, y_1,...,y_n, r)$ and dual coordinates $(\xi,r)=(\xi_1,...,\xi_{2n-1},r)$.
We have
\begin{eqnarray*}
|||u_t|||_{\epsilon}&=&|||U_t|||^2_\epsilon+\sum_{j=1}^{k}|||r_iU_t|||^2_\epsilon\ge |||U_t|||^2_\epsilon \\
&\ge&\int(1+|\xi|^2)^{\epsilon}|\hat{U_t}(x_1,...,x_{n-1}, y_1,...,y_n, r)|^2d\xi dr\\
&\ge &\int|\xi^{2\epsilon}_{2n-1}|\Big|\int\frac{\phi(y_n)\lambda(x_n)e^{-i\xi_{2n-1}y_n}dy_n}{\Big((x_n-Q(z')-1/t+iy_n\Big)^p}\Big|^2\\
&&~~~~~~~~~~\cdot\big(\Pi_{j=1}^{k}\phi(t^\epsilon_kx_j)\phi(t^{\epsilon_k}y_j)\big)^2\big(\Pi_{j=k+1}^{n-1}\phi(x_j)\phi(y_j)\big)^2dx'dx''dy'dy''d\xi_{2n-1} dr,
\end{eqnarray*}
where we use Plancherel's theorem on $\xi_1,...,\xi_{2n-2}$ in the second line. Similarly as before, we use transformations
\begin{equation*}
\begin{cases}
\tilde{x}_j=t^{\epsilon_k}x_j,\quad\tilde{y}_j=t^{\epsilon_k}y_j,~~~j=1,...,n-1,
\\
\tilde{y}_n=t y_n,\quad \tilde{\xi}_{2n-1}=1/t\xi_{2n-1}, \quad \tilde{r}=t r,
\end{cases}
\end{equation*}
and obtain
$$|||u_t|||_{\epsilon}^2\ge t^{2p-2+2\epsilon-2(n-1)\epsilon_k}J_t,$$
where 
\begin{eqnarray*}
J_t&=& \int|\tilde{\xi}_{2n-1}|^{2\epsilon}
\Big|\int\frac{\phi(t^{-1}\tilde{y}_n)\lambda(x_n(t^{-\epsilon_k}\tilde{x}_1,...,t^{-1}\tilde{r}))e^{-i\tilde{\xi}_{2n-1}\tilde{y}_n}d\tilde{y}_n}{\Big(-g+i\tilde{y}_n\Big)^p}\Big|^2\\
&&~~~~~~~~~~.\big(\Pi_{j=1}^{n-1}\phi(\tilde{x}_j)\phi(\tilde{y}_j)\big)^2d\tilde{x}'d\tilde{x}''d\tilde{y}'d\tilde{y}''d\tilde{\xi}_{2n-1} d\tilde{r}.
\end{eqnarray*}
Here $$g=-\Big(\frac{\tilde{r}-tQ(t^{-\epsilon_k}\tilde{z'})-t\sum_{j=q_0+1}^{j=n-1}|h_j(t^{-\epsilon_k}\tilde{z}_j)|^{2}}{2}-1\Big).$$
Since $tQ(t^{-\epsilon_k}\tilde z')+t\sum^{k}_{j=q_0+1}|h_j(t^{-\epsilon_k}\tilde{z}_j)|^{2}\lesssim \sum^{q_0}_{j=1}|\tilde{z}_j|^{2m}+\sum^{k}_{j=q_0+1}|\tilde{z}_j|^{2m_j}$, then if the support of $\phi$ is small enough we can assume
$$tQ(t^{-\epsilon_k}t)+t\sum^{n-1}_{j=n-k+q_0}|h_j(t^{-\epsilon_k}\tilde{z}_j)|^{2}\le 1.$$
This implies  $0<g\le \frac{-\tilde{r}+t\sum^{n-1}_{j=k+1}|h_j(t^{-\epsilon_k}\tilde{z}_j)|^{2}+3}{2} $ . 
Using a further substition 
$$y'_n=g\tilde{y}_n, ~~~~{\xi}_{2n-1}'=\frac{1}{g}\tilde{\xi}_{2n-1},$$
we get
\begin{eqnarray*}
J_t&=& \int\frac{|\xi'_{2n-1}|^{2\epsilon}}{g^{2p+1-2\epsilon}}\Big|\int\frac{\phi(\frac{y'_n}{tg})\lambda(x_n)e^{-i\xi'_{2n-1}y'_n}dy'_n}{(-1+i\tilde{y}_n)^p}\Big|^2\\
&&~~~~~~~~~~\cdot\big(\Pi_{j=1}^{n-1}\phi(\tilde{x}_j)\phi(\tilde{y}_j)\big)^2d\tilde{x}'d\tilde{x}''d\tilde{y}'d\tilde{y}''d\tilde{\xi}_{2n-1} d\tilde{r},\\
&=&J_1+J_2
\end{eqnarray*}
where $J_1$ is the integration from $-\infty$ to $-tK$, $J_2$ from $-tK$ to $0$ and where  $K$ is suitably chosen. Note that $J_1\ge 0$. Now, we consider $J_2$.\\

For $\tilde{r}\in [-tK, 0]$, we see that
$$|x_n|=|\frac{\tilde{r}/t+Q(t^{-\epsilon_k}\tilde{z}')-\sum^{n-1}_{j=q_0+1}|h_j(t^{-\epsilon_k}\tilde{z}_j)|^2}{2}|\le C.$$

We may choose $\lambda\in C_0(\mathbb{R})$ such that $\lambda(x)=1$ for $|x|\le C$. Then
$$\int|\xi'_{2n-1}|^{2\epsilon}\Big|\int\frac{\phi(\frac{y'_n}{tg})\lambda(x_n)e^{-i\xi'_{2n-1}y'_n}dy'_n}{(-1+i\tilde{y}_n)^p}\Big|^2d\xi'_{2n-1}\ge \text{const}>0.$$
It follows
\begin{eqnarray*}
J_2&\gtrsim&\int\int^0_{\tilde{r}=-tK}\frac{\big(\Pi_{j=1}^{n-1}\phi(\tilde{x}_j)\phi(\tilde{y}_j)\big)^2d\tilde{x}'d\tilde{x}''d\tilde{y}'d\tilde{y}'' } {\Big(-\tilde{r}+t\sum^{n-1}_{j=k+1}|h_j(t^{-\epsilon_k}\tilde{z}_j)|^{2}+3\Big)^{2p+1-2\epsilon}}d\tilde{r}\\
&\gtrsim&\int\frac{\big(\Pi_{j=1}^{n-1}\phi(\tilde{x}_j)\phi(\tilde{y}_j)\big)^2d\tilde{x}'d\tilde{x}''d\tilde{y}'d\tilde{y}'' } {\Big(t\sum^{n-1}_{j=k+1}|h_j(t^{-\epsilon_k}\tilde{z}_j)|^{2}+3\Big)^{2p-2\epsilon}}d\tilde{r}\\
&&~~~~~~~~~~-\int\frac{\big(\Pi_{j=1}^{n-1}\phi(\tilde{x}_j)\phi(\tilde{y}_j)\big)^2d\tilde{x}'d\tilde{x}''d\tilde{y}'d\tilde{y}'' } {\Big(tK+t\sum^{n-1}_{j=k+1}|h_j(t^{-\epsilon_k}\tilde{z}_j)|^{2}+3\Big)^{2p-2\epsilon}}d\tilde{r}\\
&\gtrsim&\int\frac{\big(\Pi_{j=1}^{n-1}\phi(\tilde{x}_j)\phi(\tilde{y}_j)\big)^2d\tilde{x}'d\tilde{x}''d\tilde{y}'d\tilde{y}'' } {\Big(t\sum^{n-1}_{j=k+1}|h_j(t^{-\epsilon_k}\tilde{z}_j)|^{2}+3\Big)^{2p-2\epsilon}}d\tilde{r}.\\
\end{eqnarray*}
The last inequality follows from the fact that we can choose $K$ and $t$ such that
\begin{eqnarray*}
tK+t\sum^{n-1}_{j=k+1}|h_j(t^{-\epsilon_k}\tilde{z}_j)|^{2}+3\ge 2(t\sum^{n-1}_{j=k+1}|h_j(t^{-\epsilon_k}\tilde{z}_j)|^{2}+3).
\end{eqnarray*}

Then 
\begin{eqnarray*}
J_t\gtrsim\int_0^{\delta}...\int_0^\delta\frac{d\tilde{x}_{k+1}d\tilde{y}_{k+1}...d\tilde{x}_{n-1}d\tilde{y}_{n-1}}{(t\sum_{j=k+1}^{n-1}|t^{-\epsilon_k}z_j|^{2m_j}+1)^{2p-2\epsilon}}\cong t^{-\sum_{j=k+1}^{n-1}\frac{1}{m_j}+2(n-k-1)\epsilon_k},\\
\end{eqnarray*}
where the last inequality follows by Lemma~\ref{l6.1}.
So we have 
\begin{equation}
\Label{new}
|||u_t|||^2_\epsilon  \gtrsim t^{2p-2+2\epsilon-2k\epsilon_k-\sum_{j=k+1}^{n-1}\frac{1}{m_j}}.
\end{equation}
Since subelliptic estimates hold with order $\epsilon$ for any $k$-form ($q_0+1\le k\le n-1$), then
\begin{equation}
\Label{6.10}
|||u_t|||_\epsilon^2\lesssim Q(u_t,u_t).
\end{equation}
Combining \eqref{6.9} \eqref{new} and \eqref{6.10}, we get $\epsilon \le \epsilon_k$

The proof of Theorem~\ref{t1.3} (i) is complete.

\hskip14cm$\Box$

\noindent
{\it Proof of Theorem~\ref{t1.3} (ii)}
We proceed in similar way as in the proof of Theorem~\ref{t1.3} (i) and choose the coefficient of our form by setting 
$$\begin{cases}
f_t(z',z_n)=(z_n-\sum_{j=1}^{q_0}|h_j(z_j)|^2-1/t)^{-p}\\
\Phi_t(z)=\big(\Pi_{j=1}^{n-1}\phi(t^{\epsilon_k}x_j)\phi(t^{\epsilon_k}y_j)\big)\lambda(x_n)\phi(y_n)
\end{cases}$$

Then 
\begin{equation*}
\begin{split}
Q(u_t,u_t)&\lesssim \sum_{j=1}^{q_o}\no{r_{z_j}\frac{\partial f_t}{\partial \bar{z}_j}\Phi_t}^2+\sum_{j=k+1}^{q_o}\no{\frac{\partial f_t}{\partial \bar{z}_j}\Phi_t}^2
\\
&+\underset{j=1,...,k,i}{\underset{{i=k+1,...,n-1}}\sum}|||r_{z_j}||r_{z_i}|\partial_{z_n}(f_t)\Phi_t||^2+\sum_{j=1}^{n}\no{f_t\frac{\partial \Phi_t}{\partial z_j}}^2+\no{U_t}^2.
\end{split}
\end{equation*}
We can show that
$Q(u_t,u_t)\lesssim t^{2p-2+2\epsilon_k-2(n-1)\epsilon_k}I_t$ where 
$$
I_t=\int_0^\delta...\int_0^\delta\frac{dx_1dy_1...dx_kdy_k}{\left(t\underset {j=1}{\overset k \sum}|t^{-\epsilon_k}z_j|^{2m_j}+1\right)^{2p-2}}.
$$
Owing to Lemma~\ref{l6.1} we have $I_t\simle t^{-\sum_{j=1}^k \frac 1{m_j}+2k\epsilon_k}$ which yields
$$Q(u_t,u_t)\lesssim t^{2p-2+2\epsilon_k-2(n-k-1)\epsilon_k-\sum_{j=1}^{k}\frac{1}{m_j}}.$$
Similarly, we have
$$|||u_t|||_\epsilon^2\gtrsim t^{2p-2+2\epsilon_k-2(n-k-1)\epsilon_k-\sum_{j=1}^{k}\frac{1}{m_j}},$$
which yields the conclusion of the proof of Theorem~\ref{t1.3} (ii).

\hskip14cm$\Box$

\end{document}